\documentclass[final,3p,times,authoryear]{elsarticle}
\usepackage{amssymb}
\usepackage{color,latexsym,amsfonts,amssymb}
\usepackage{amsmath}
\usepackage{amsfonts,amssymb}
\usepackage{epsfig,epstopdf}
\usepackage{bm}
\usepackage{amsfonts,amsmath,latexsym}
\usepackage{mathrsfs}
\usepackage{amssymb,amsmath,bm}

\usepackage{ulem}%下划线

\usepackage{geometry}
\usepackage{caption}
\usepackage{enumerate}
\usepackage{url}
\usepackage{graphicx}
\usepackage{subfigure}

\usepackage{booktabs}
\usepackage{threeparttable}
\usepackage{multirow}
\usepackage{algorithm}
\usepackage{algorithmic}
\usepackage{verbatim}% for begin{comment}
\newtheorem{lemma}{Lemma}[section]
\newtheorem{definition}{Definition}[section]

\newtheorem{theorem}{Theorem}[section]
\newtheorem{assumption}{Assumption}[section]

\newtheorem{condition}{Condition}[section]

\def\proclaim#1{\par \bigskip\noindent {\bf #1}\bgroup\it\ }
\def\endproclaim{\egroup\par\bigskip}
\def\proof{\par\noindent{\bf Proof.} \;}

\newbox\TempBox \newbox\TempBoxA

\def\ep{\mathbb{E}}

\def\R{\mathbb{R}}
\def\I{\mathbb{I}}
\def\X{\mathbf{X}}
\def\x{\mathbf{x}}
\def\z{\mathbf{z}}
\def\P{\mathbb{P}}
\def\bbeta{\boldsymbol{\beta}}

\def\bomega{\boldsymbol{\omega}}

\newcommand{\bq}{\begin{equation}}
\newcommand{\eq}{\end{equation}}
\newcommand{\bqs}{\begin{equation*}}
\newcommand{\eqs}{\end{equation*}}
\newcommand{\be}{\begin{eqnarray}}
\newcommand{\ee}{\end{eqnarray}}
\newcommand{\by}{\begin{eqnarray*}}
\newcommand{\ey}{\end{eqnarray*}}
\newcommand{\bn}{\begin{enumerate}}
\newcommand{\en}{\end{enumerate}}
\newcommand{\bi}{\begin{itemize}}
\newcommand{\ei}{\end{itemize}}
\newcommand{\bds}{\begin{description}}
\newcommand{\eds}{\end{description}}
\newcommand{\bcen}{\begin{center}}
\newcommand{\ecen}{\end{center}}
\journal{}

\begin{document}

\begin{frontmatter}



\title{Robust estimation and shrinkage in ultrahigh dimensional expectile regression with heavy tails and variance heterogeneity\tnoteref{label_title}}
\tnotetext[label_title]{This research is partly supported by the Fundamental Research Funds for the Central Universities, Major Project of the National Social Science Foundation of China (No.13$\&$ZD163) , Zhejiang Provincial Natural Science Foundation (No: LY18A010005) and the Research Project of Humanities and Social Science of Ministry of Education of China(No. 17YJA910003).}

\author[Label1]{Jun Zhao}
%\fntext[Label1]{Zhejiang University City College}

\author[Label2]{Guan'ao Yan}

\author{Yi Zhang\corref{cor1}\fnref{Label2}}
%\fntext[Label1]{Department of Statistics, Zhejiang University City College}
%\ead{ }
\address[Label1]{Zhejiang University City College}
\address[Label2]{School of Mathematical Sciences, Zhejiang University}
\cortext[cor1]{Address for correspondence: Yi Zhang, School of Mathematical Sciences, Zhejiang University, 38 Zheda Road, Hangzhou, 310027, Zhejiang, China.}
\ead{zhangyi63@zju.edu.cn}

\begin{abstract}
High-dimensional data subject to heavy-tailed phenomena and heterogeneity are commonly encountered in various scientific fields and bring new challenges to the classical statistical methods. In this paper, we combine the asymmetric square loss and huber-type robust technique to develop the robust expectile regression for ultrahigh dimensional heavy-tailed heterogeneous data. Different from the classical huber method, we introduce two different tuning parameters on both sides to account for possibly asymmetry and allow them to diverge to reduce bias induced by the robust approximation. In the regularized framework,
we adopt the generally folded concave penalty function like the SCAD or MCP penalty for the seek of bias reduction. We investigate the finite sample property of the corresponding estimator and figure out how our method plays its role to trades off the estimation accuracy against the heavy-tailed distribution. Also, noting that the robust asymmetric loss function is everywhere differentiable, based on our theoretical study, we propose an efficient first-order optimization algorithm after locally linear approximation of the non-convex problem. Simulation studies under various distributions demonstrates the satisfactory performances of our method in coefficient estimation, model selection and heterogeneity detection.

\end{abstract}


\end{frontmatter}


\section{Introduction}
\label{intro}

%% outline for the introduction part
%% 1.1 Expectile的用处及最近的研究进展，从而说明Expectile的优势；
%% 1.2 Expectile在处理重尾情形时存在的问题，即 维度问题 $p=O(n^\alpha)$ （回归误差的矩条件）以及 Error Bound （模拟和理论上说明）受误差矩条件的影响；
%% 1.3 Robust 的引入，包括形式以及其他文献中类似的处理；
%% 1.4 文章的贡献
%% 1.5 文章的结构安排
%% 异质性的处理放在哪里？
Inspired by the asymmetric check loss in quantile regression, \cite{AAP76} and \cite{New87} used the asymmetric square loss instead in linear regression and proposed the so-called expectile regression.
 Analogous to quantile regression, expectile regression assigns different weights onto positive and negative squared error losses   respectively so that it can draw a complete picture
of the conditional distribution of the response variable given the covariates, making it a useful tool for modeling heterogeneous data.   However, expectile regression does have certain advantages over quantile regression in both theoretical and computational aspects. First, its loss function $\phi_{\alpha}(\cdot)$ is everywhere differentiable so that  algorithms based on  first-order optimization condition can be applied to alleviate the computation burden, especially in the high dimensional setting. Second, due to differentiability of $\phi_{\alpha}(\cdot)$, estimation in expectile regression is straightforward and  the asymptotic covariance matrix of the expectile regression
estimator does not involve estimating the density function of the errors (\cite{New87}), which may bring convenience when the density function is hard to estimate.  Other comparisons between quantile and expectile regression can be found in \cite{Wal15}.
Because of these features, expectile regression has made considerable development in the last two decades, for example, \cite{Efr91}, \cite{YaT96},\cite{DRH09},\cite{SaK12},\cite{KaL16} and among others.

% Signal to noise ratio can be introduced to explain our motivation.....

It has to be said that compared to quantile regression, while expectile regression has those advantages above, it still has some shortcomings.
Expectile regression is more sensitive to the extreme values in the response variable than quantile regression, which is more obvious in high dimension. At the same time, nowadays, massive data subject to heavy-tailed distributions are frequently observed in many scientific areas, from microarray experiments to finance (\cite{FLW17}) and bring new challenges to conventional statistical methods. \cite{ZCZ18} apply expectile regression for analyzing heteroscedasticity in high-dimensional linear model and point out that  parameter estimation and variable selection of expectile regression in high dimension are affected by the regression error moment order. Moreover, the dimensionality expectile regression can handle can only increase with some certain polynomial rate of the sample size, also affected by the error moment order.  Huber-type robust  methods provide potential solutions to tackle the problem of heavy-tailed errors in this situation. \cite{FLW17} and \cite{ZBFL18} proposed Huber-type
estimators in both low and high dimensional linear mean regression and derived non-asymptotic deviation bounds for the estimation error. However, no literature has focused on expectile regression for ultrahigh dimensional heavy-tailed and heterogeneous data.

In this paper, we combine the asymmetric square loss and huber-type robust technique to develop the robust expectile regression for ultrahigh dimensional heavy-tailed heterogeneous data. Different from the classical huber estimation, we introduce two different tuning parameters on both sides to account for possibly asymmetric error distribution and allow them to diverge to reduce bias induced by the robust approximation. In the regularized framework,
we adopt the generally folded concave penalty function like the SCAD or MCP penalty for the seek of bias reduction. We investigate the finite sample property of the corresponding estimator and figure out  how  our method plays the trick to trades off the estimation accuracy against the heavy-tailed distribution with a little cost. Also, noting that the robust asymmetric loss function is everywhere differentiable, we propose an efficient algorithm based on first-order optimization condition after locally linear approximation of the nonconvex penalty. Various simulations demonstrates the good performances of our method.

The remainder of this article is organized as follows. In Section 2, we introduce the robust asymmetric square loss and the generally folded concave penalty function, and on this basis, we further propose the adjusted regularized expectile regression. An efficient algorithm for the corresponding optimization problem is also provided.
In Section 3, we present our main theoretical results, including the finite sample property of out estimator.
% and also state how  our method plays the trick to  trades off  the estimation error against the heavy-tailed distribution with a little cost.
In Section 4, we conduct simulation studies under various shapes of error distribution.  Technical proofs are presented in the Appendix.


\section{Methodology}
\label{method}
\setcounter{definition}{0}\setcounter{definition}{0}
\setcounter{equation}{0}\setcounter{lemma}{0}
\setcounter{proposition}{0}\setcounter{theorem}{0}
\setcounter{remark}{0}\setcounter{corollary}{0}
\subsection{ Robust Asymmetric Loss}
As given out in \cite{AAP76} and \cite{New87},
the asymmetric square loss function $\phi_{\alpha}(\cdot)$ is defined as follows,
\be\label{phi}
\phi_{\alpha}(r)=|\alpha-\I(r<0)|r^2=
\begin{cases}
	\alpha r^2, &r \geq 0, \\
	(1-\alpha)r^2, &r<0,
\end{cases}
\ee
where $\alpha$ controls the degree of loss asymmetry and is called the expectile level.
Denote the $\alpha$-th expectile of random variable $y$  by
$ m_{\alpha}(y)=\underset{m\in\R}{\arg\min}~ \ep\phi_{\alpha}(y-m)$.
Note that 1/2-th expectile is exactly the mean.

%??\cite{New87} point out that expectile regression is a useful tool for modelling heterogeneous data, especially for testing homoscedasticity hypothesis. In high-dimensional situations, \cite{ZCZ18} applied expectile regression for analyzing heteroscedasticity in high-dimensional linear model.

\begin{figure}[!htb]
	\centering
	\subfigure[Different Loss Functions at $\alpha=0.25$]{
	\includegraphics[width=0.45\textwidth]{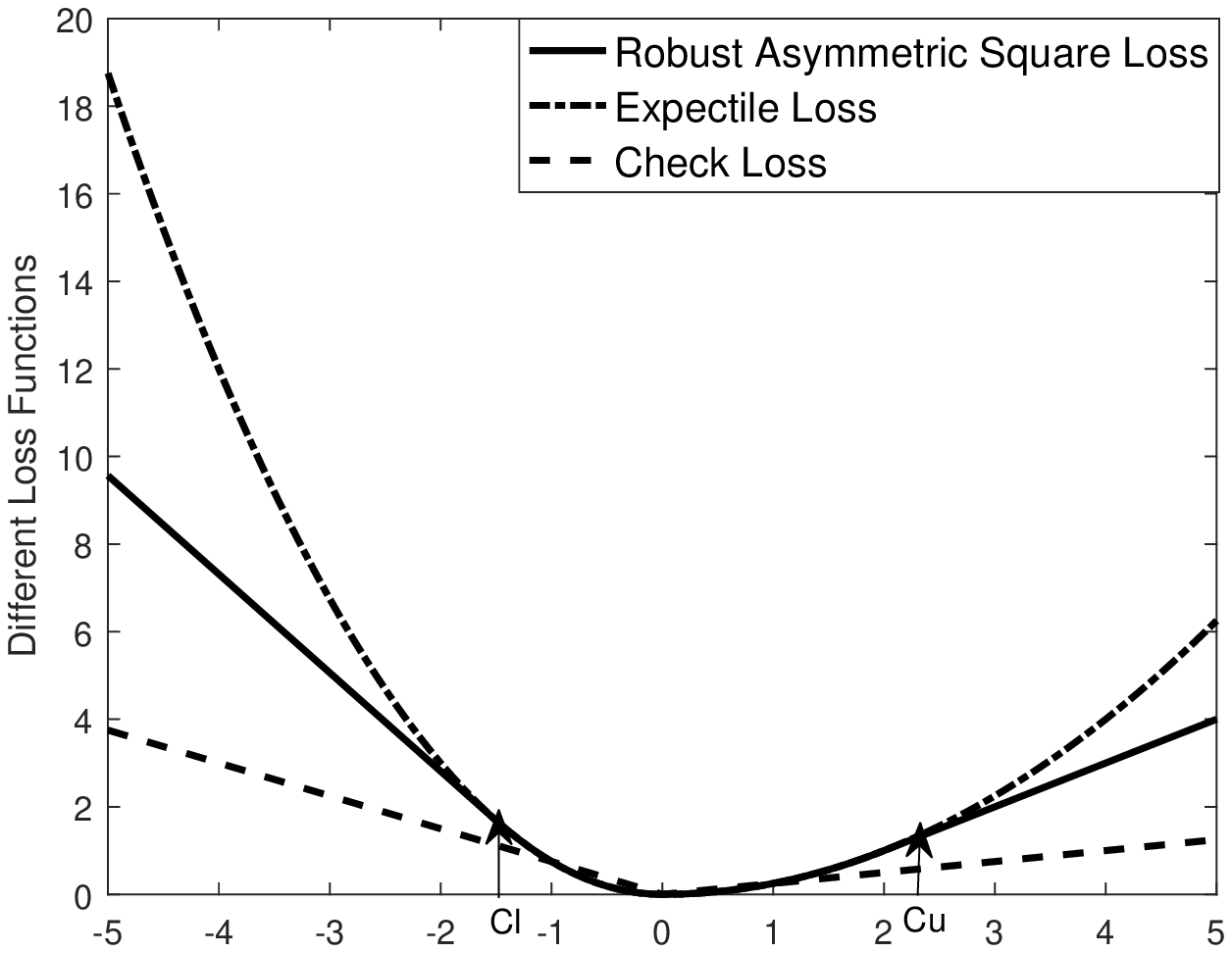}}
	\subfigure[Different Loss Functions at $\alpha=0.75$]{
	\includegraphics[width=0.45\textwidth]{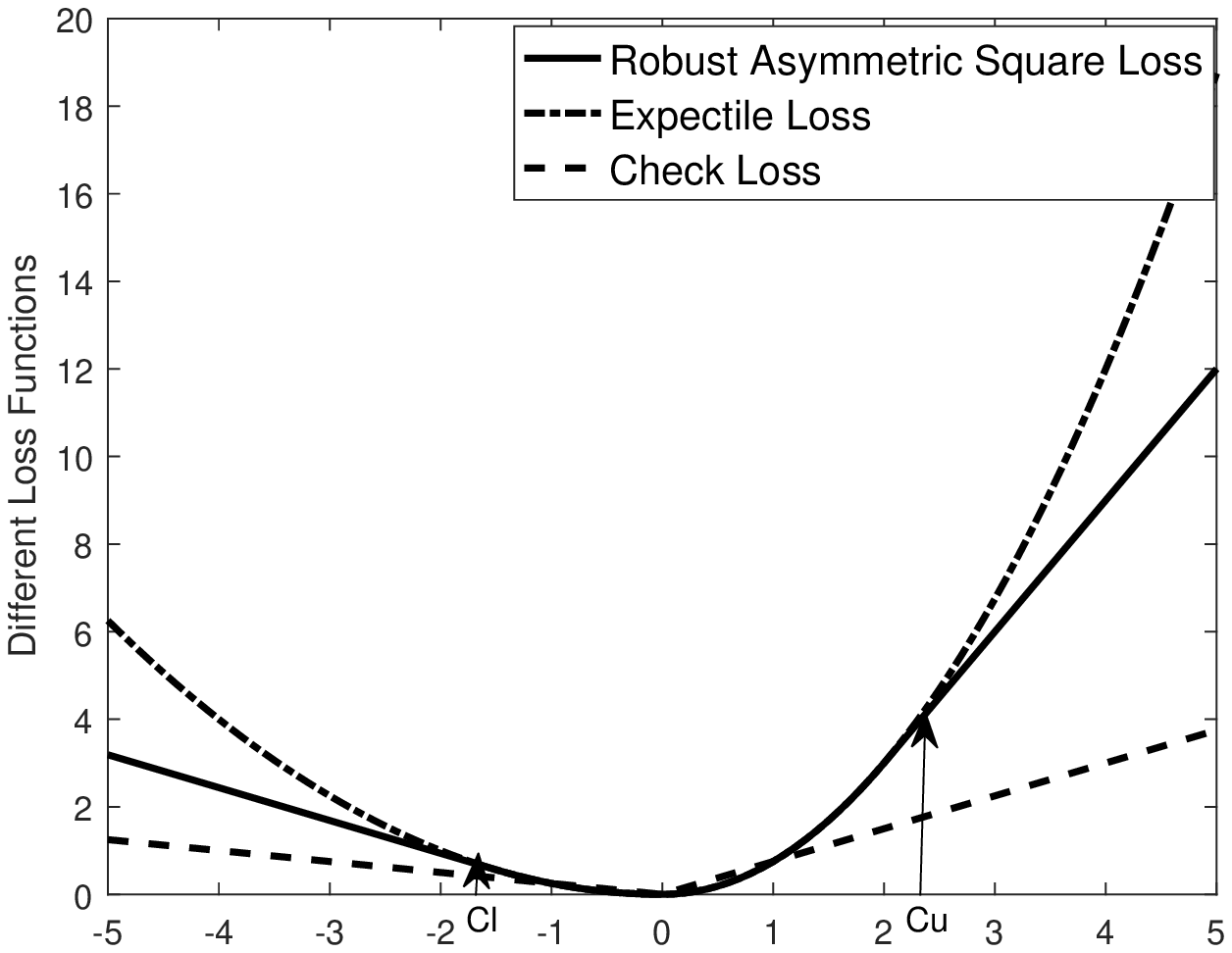}}
	\caption{Different Loss functions }\label{fig:Robust Asymmetric Loss}
\end{figure}

To tackle the problem of heavy-tailed errors, huber-type robust technique provides potential solutions.
The classical huber loss (\cite{Hub64}) is a hybrid of squared loss for relatively small errors and absolute
loss for relatively large errors, where the degree of hybridization is controlled by one tuning
parameter. In our paper, considering the role of expectile regression in analyzing heterogeneity, we define the following robust asymmetric loss function with two different tuning parameters,
\be\label{huber phi}
\psi_{\alpha}(r;C_u,C_l)=
\begin{cases}
	2\alpha C_u r-\alpha C_u^2, &r \geq C_u, \\
	\alpha r^2, & 0\leq r<C_u,\\
	(1-\alpha)r^2, &C_l<r<0,\\
	2(1-\alpha)C_lr-(1-\alpha)C_l^2, & r\leq C_l,
\end{cases}
\ee
See Figure \ref{fig:Robust Asymmetric Loss}.
This loss function is quadratic for small values of $r$ and linear for large values of $r$, sharing the same robust idea with the classical huber loss. What a significant difference between these two loss functions is there exist two different truncation tuning parameters $C_u$ and $C_l$ in the robust asymmetric loss, in accordance with "asymmetry" and can diverge for bias reduction. Expectile asymmetric square loss and quantile check loss can be regarded as the two particular scenarios when $C_u$ and $C_l$ both go to $\infty$ or 0 in some sense. The truncation parameters $C_u$ and $C_l$ control the trade-off balance between bias and robustness and need to be determined by some data-driven methods.

\subsection{Regularized Framework}
 	
Let us begin with the notation and statistical setup. Suppose that we have a high-dimensional data sample $\{Y_i,\x_i\}$,~$i=1,\ldots,n$, where $\x_i=(x_{i1},\ldots,x_{ip}),~i=1,\ldots,n$ are independent and identically distributed $p$-dimensional covariates along with the common mean 0. Consider the data are generated from the following linear model,
\bq\label{high-dimensional linear model}
Y_i=\beta_1^*x_{i1}+\cdots+\beta_p^*x_{ip}+\epsilon_i=\x_i'\boldsymbol{\beta}^*+\epsilon_i.
\eq

To account for data heterogeneity, here we introduce the so-called "variance heterogeneity" into model (\ref{high-dimensional linear model}) from \cite{RaS96} in their "mean and dispersion additive model". Specifically,  $\epsilon_i$ can take the following form, $$\epsilon_i=\sigma(\x_i)\eta_i,$$ where $\sigma(\x_i)$, conditional on $\x$, can be nonconstant, linear form (\cite{Ait87},~\cite{GaZ16}), nonparametric form (\cite{RaS96}).

In addition to heterogeneity, $\{\epsilon_i\}_{i=1}^n$ are  assumed to be mutually independent and satisfy $m_{\alpha}(\epsilon_i|\x_i)=0$ for some specific $\alpha$.
Thus, the $\alpha$-th conditional expectile of $Y_i$ given covariates $\x_i$ through the model (\ref{high-dimensional linear model}) is
 $m_{\alpha}(Y_i|\x_i)=\x_i'\boldsymbol{\beta}^*$.
So $\bbeta^*$  actually minimizes the conditional population risk $\ep[\phi_{\alpha}(Y_i-\x'\bbeta)|\x]$, i.e.,
\be
\bbeta^*=\underset{\bbeta\in\R^p}{\arg\min} \ep[\phi_{\alpha}(Y_i-\x'\bbeta)|\x],
\ee
and is assumed to be unique for convenience. In consideration of data heterogeneity, we should address that $\bbeta^*$ may change for different weight levels.

% How to determine the satisfied $\alpha$ can be refered to  Yao and Tong 1996 Proposition 1????????

Here we consider the scenario where $p=p(n)$ increases with the sample size $n$ at the exponential rate $\log (p)=O(n^b)$ for $0<b<1$,  i.e. the ultra-high dimensional settings. One leading way to deal with this situation is to assume that the true parameter $\boldsymbol{\beta}^{*}=(\beta_1^{*},\ldots,\beta_p^{*})$ is sparse. Let $A=\{j: \beta_j^{*}\neq0,1\leq j\leq p\}$ be the active index set and its cardinality $q=q(n)=|{A}|$. Sparsity means that  $q<n$ and all the left $(p-q)$ coefficients are exactly zero. Generally, similar to $\boldsymbol{\beta}^{*}$, $A$ and $|{A}|$ may also change for different expectile levels and for simplicity in notation, we omit such dependence when no confusion arises. Without loss of generality, we rewrite $\boldsymbol{\beta}^*=((\boldsymbol{\beta}_A^*)',\mathbf{0}')'$ where $\boldsymbol{\beta}_A^*\in \R^q$ and $\mathbf{0}$ denotes a $(p-q)$ dimensional vector of zero. Let $\X=(\x_1,\ldots,\x_n)'$ be the $n\times p$ matrix of covariates. Denote $\X_j$ the $j$-th column of $\X$ and define $\X_A$ the submatrix of $\X$ that consists of its first $q$ columns.
% and denote by $\z_i$ the $i$-th row of $\X_A$.

Under sparsity assumption,
regularized framework has been playing a leading role in analyzing high-dimensional data in the past two decades.
Throughout this paper, we assume that the regularizer $P_{\lambda}(t)$ is a generally
folded concave penalty function defined on $t \in (-\infty,\infty)$ satisfying:
\begin{assumption}\label{assumption on penalty}
\bn[(i)]
\item $P_{\lambda}(t)$ is symmetric around the origin, i.e., $P_{\lambda}(t)=P_{\lambda}(-t)$ with $P_{\lambda}(0)=0$;
\item $P_{\lambda}(t)$ is nondecreasing  in $t\in (0,\infty)$, but the function $\frac{P_{\lambda}(t)}{t}$ is nonincreasing $t\in (0,\infty)$;
\item $P_{\lambda}(t)$  is differentiable $\forall~ t\neq 0$ and subdifferentiable at $t=0$ with $\lim_{t\rightarrow0^+}P'_{\lambda}(t):= \lambda L$;
%\item  $P'_{\lambda}(t)\geq a_1\lambda$ for $t\in (0,a_2\lambda)$ and $P'_{\lambda}(0):= a_1\lambda$ ;
\item  There exists $\mu>0$ such that $P_{\lambda,\mu}(t):=P_{\lambda}(t)+\frac{u}{2}t^2$ is convex.
\en
\end{assumption}
These assumptions about the regularizer are commonly used in nonconvex regularized framework, for example, see \cite{LaW15}.
Many popular non-convex regularizers in practice satisfy all the assumptions above, like the SCAD or MCP penalty, details in Appendix A.2 in \cite{LaW15}.
\bi
\item \textbf{SCAD, \cite{Fan01}}. ~The SCAD penalty is defined through its first order derivative and symmetry around the origin. To be specific, for $\theta > 0$,
\be
P_{\lambda}'(\theta)=\lambda\{\I(\theta\leq\lambda)+\frac{(a\lambda-\theta)_+}{(a-1)\lambda}\I(\theta>\lambda)\},
\ee
where $a> 2$ is a fixed parameter. By straight calculation, for $\theta > 0$,
 \be
 P_{\lambda}(\theta)&=&\lambda\theta\I(\theta\leq\lambda)+\frac{a\lambda\theta-(\theta^2+\lambda^2)/2}{a-1}\I(\lambda\leq\theta\leq a\lambda)+\frac{(a+1)\lambda^2}{2}\I(\theta> a\lambda).
 \ee
 The SCAD penalty satisfies all the conditions with $L=1$ and $\mu=\frac{1}{a-1}$.
\item \textbf{ MCP, \cite{Zha10}}.~The MCP penalty has the following form:
\be
P_{\lambda}(\theta)=\text{sgn}(\theta)\lambda\int_{0}^{|\theta|}(1-\frac{z}{\lambda b})dz,
\ee
where $b>0$ is a fixed parameter and $\text{sgn}(\cdot)$ is the sign function. The MCP penalty satisfies all the conditions with $L=1$ and $\mu=\frac{1}{b}$.
\ei

  Finally, we define the following robust regularized expectile loss with the non-convex penalty:
\bq\label{penalty loss function}
L(\boldsymbol{\beta})=\frac{1}{n}\sum_{i=1}^n\psi_{\alpha}(y_i-\x_i'\boldsymbol{\beta};C_u,C_l)+\sum_{j=1}^pP_{\lambda}(|\beta_j|),
\eq
where $\psi_{\alpha}(\cdot;C_u,C_l)$ is defined by (\ref{huber phi}).
Then the estimate of the coefficients is obtained by solving the optimization problem below,
\bq\label{non-convex regularized problem}
\hat{\boldsymbol{\beta}}=\underset{\boldsymbol{\beta} \in \R^{p}}{\arg\min}~L(\boldsymbol{\beta}).
\eq

\subsection{Optimization Algorithm}
Due to the non-convexity of the penalty, optimization problem (\ref{non-convex regularized problem}) is a non-convex  one in high dimension. Here we take use of the Local Linear Approximation (LLA,~\cite{ZaL08}) strategy, due to its computational efficiency and good statistical properties (see \cite{FXZ14}), to approximate the penalized optimization problem into a convex one. Different from one-step optimization procedure in \cite{ZaL08}, in our paper, we use the iteration strategy to seek for the "best" minimizer.
The algorithm uses the adaptive weight $\bomega^{(t)}$ in its every iteration and solves a sequence of convex programs.
The details are summarized in Algorithm \ref{LLA algorithm}.
\begin{algorithm}[!htb]
\caption{The adaptive LLA algorithm for the non-convex optimization problem (\ref{non-convex regularized problem})}
\label{LLA algorithm}
\begin{algorithmic}[1]
\STATE Initialize $\boldsymbol{\beta}=\boldsymbol{\beta}^{(0)}$.
\STATE For $t=1,2,\ldots$, repeat the following iteration (a),(b) and (c) until convergence
%($\sum_{j=1}^p|\boldsymbol{\beta}^{t+1}-\boldsymbol{\beta}^t|$ is sufficiently small)
\bn[(a)]
\item At step $t$, calculate the corresponding weight based on the current solution $\bbeta^{(t-1)}=(\beta_1^{(t-1)},\ldots,\beta_p^{(t-1)})'$
\bqs \bomega^{(t)}=(\omega_1^{(t)},\ldots,\omega_p^{(t)})'=(P_{\lambda}'(|\beta_1^{(t-1)}|),\ldots,P_{\lambda}'(|\beta_p^{(t-1)}|))'.\eqs
\item So, at step $t$, calculate the current local linear approximation of  regularized expectile loss function $L(\boldsymbol{\beta})$, denoted by $L(\boldsymbol{\beta}|\boldsymbol{\beta}^{(t-1)})$,
\bqs
L(\boldsymbol{\beta}|\boldsymbol{\beta}^{(t-1)})=\frac{1}{n}\sum_{i=1}^n\psi_{\alpha}(y_i-\x_i'\boldsymbol{\beta};C_u,C_l)+\sum_{j=1}^p\omega_j^{(t)}|\beta_j|.
\eqs
\item $\boldsymbol{\beta}^{(t)}=\underset{\boldsymbol{\beta}\in \R^{p}}{\arg\min}~ L(\boldsymbol{\beta}|\boldsymbol{\beta}^{(t-1)})$.

\en
\end{algorithmic}
\end{algorithm}

There are many choices for the initial value $\boldsymbol{\beta}^{(0)}$. For example, $\boldsymbol{\beta}^{(0)}$ can be chosen as the estimator obtained from the convex Lasso-type penalized expectile regression to accelerate the convergence of our algorithm. Or simply and roughly, the worst choice $\bbeta^{(0)}=\boldsymbol{0}$ is also available. In fact, if $\bbeta^{(0)}=\boldsymbol{0}$, then in the following, $\bomega^{(1)}=(P_{\lambda}'(0^+),\ldots,P_{\lambda}'(0^+))'=(\lambda L,\ldots,\lambda L)'$, which also results in the Lasso-type framework.
The robust asymmetric square loss function $\psi_{\alpha}(r;C_u,C_l)$ is everywhere differentiable, resulting in convenience in computation. So, algorithms based on the first-order optimization condition, for example, proximal gradient method (\cite{PaB14}), can be applied to solve the corresponding problem (c) in \textbf{Algorithm \ref{LLA algorithm}}. And for such kind of convex optimization problem with differentiability in \textbf{Algorithm \ref{LLA algorithm}}, there are many available and powerful programs to solve it. For example, in our simulation, we use CVX, a Matlab optimization tool for specifying and solving convex programs; see \cite{MaS08},~\cite{MaS13}.
\begin{comment}
\subsection{Algorithm Convergence}
Robust oracle estimate is defined through the following optimization problem:
\be\label{Oracle Estimator Representation}
\hat{\bbeta}^*=\underset{\bbeta}{\arg\min}\frac{1}{n}\sum_{i=1}^n\psi_{\alpha}(y_i-\z_i'\bbeta;C_u,C_l).
\ee
\begin{theorem}
The strong oracle property of the LLA algorithm, i.e., the proposed LLA algorithm converges to the oracle estimator after a few steps with probability tending to 1.
\end{theorem}
\textbf{Paper I}:
\bi
\item Question 1: The oracle estimator is one of the local minimizers which satisfies the first order optimization necessary condition;
    \item Question 2: The theoretical estimation bound for the stationary points in regularized nonconvex robust expectile regression.
    \item Question 3: The composite gradient methods in Matlab, here whether to determine the convergence rate remains discussion;
\ei
\textbf{Paper II}:
\bi
\item Question 1: The error bound between the oracle estimate and the theoretical optima by the regularized oracle expectile regression;

\item Question 2: Convergent rate for lasso-type estimators since the LLA algorithm solves weighted lasso-type optimization problems at every iteration; see \cite{ANW12}

\item Question 3: The probability at which the LLA algorithm converges to the oracle estimator;
\ei
\end{comment}

\section{Main Results}
\setcounter{definition}{0}\setcounter{definition}{0}
\setcounter{equation}{0}\setcounter{lemma}{0}
\setcounter{proposition}{0}\setcounter{theorem}{0}
\setcounter{remark}{0}\setcounter{corollary}{0}

The true coefficients $\bbeta^*$  minimizes the population risk $\bbeta^*=\underset{\boldsymbol{\beta} \in \R^{p}}{\arg\min}~\ep\phi_{\alpha}(y-\x'\bbeta)$. If  Condition \ref{C(2)} below holds for example, then $\ep\phi_{\alpha}(y-\x'\bbeta)$ is strongly convex. By convex analysis, $\bbeta^*$ actually satisfies the sufficient and necessary condition,
\bq\label{theoretical optimal solution condition}
\ep \nabla \phi_{\alpha}(y-\x'\bbeta^*)=0.
\eq
and therefore $\bbeta^*$ is globally unique. For non-convex optimization problem (\ref{non-convex regularized problem}), local minima $\hat{\boldsymbol{\beta}}$ ( may not be unique) must satisfy the usual first-order necessary subgradient condition, i.e., for $\hat{\beta_j}$
\be\label{stationary solution condition}
-\frac{1}{n}\sum_{i=1}^n\psi'_{\alpha}(y_i-\x_i'\hat{\bbeta};C_u,C_l)\x_{ij} + P_{\lambda}'(\hat{|\beta_j|})=0,
\ee
where $\psi'_{\alpha}(\cdot)$ is the first derivative of the robust asymmetric square loss function,
\be\label{huber phi derivative}
\psi'_{\alpha}(r;C_u,C_l)=
\begin{cases}
2\alpha C_u, &r \geq C_u, \\
2\alpha r, & 0\leq r<C_u,\\
2(1-\alpha)r, &C_l<r<0,\\
2(1-\alpha)C_l, & r\leq C_l.
\end{cases}
\ee

In this section, we establish the finite sample error bound for any estimates satisfying  equation (\ref{stationary solution condition}).
Since we introduce the robust technique in expectile regression, denote by \textbf{RER} (\textbf{R}obust \textbf{E}xpectile \textbf{R}egression) for short, to handle the scenario where heavy-tailed phenomena possibly exists, we define the approximate pseudo true coefficients
\be
\bbeta^*(C_u,C_l)=\underset{\boldsymbol{\beta} \in \R^{p}}{\arg\min}~\ep\psi_{\alpha}(y-\x'\bbeta;C_u,C_l).
\ee
Then the statistical error $\|\hat{\bbeta}-\bbeta^*\|_2$ can be bounded by
\bq
\|\hat{\bbeta}-\bbeta^*\|_2\leq \|\bbeta^*(C_u,C_l)-\bbeta^*\|_2+\|\hat{\bbeta}-\bbeta^*(C_u,C_l)\|_2,
\eq
where we call the first and second term on the right side the approximation error and estimation error respectively.

The following theorem gives out the upper bound of the approximation error.
Before we state our theoretical results, we list some needed conditions.
\begin{condition}\label{C(1)}
For the regression error,  $\ep[\epsilon^k|\x]\leq M_k<\infty$ for some $k>2$.
\end{condition}
\begin{condition}\label{C(2)}
 There exist $\kappa_l,\kappa_u>0$ for $\Sigma:=\ep \x\x'$ such that $0<\kappa_l\leq \lambda_{\min}(\Sigma)\leq \lambda_{\max}(\Sigma)\leq \kappa_u<\infty$ where $\lambda_{\min}(\cdot)$ and $\lambda_{\max}(\cdot)$ represent the eigenvalue operator
\end{condition}
\begin{condition}\label{C(3)}
For any $\nu\in \R^p $, $\x'\nu$ is sub-Gaussian distributed with parameter at most $\kappa_0^2\|\nu\|_2^2$, i.e., $\ep{\exp(t\x'\nu)}\leq \exp(t^2)\kappa_0^2\|\nu\|_2^2/2$, for any $t\in R^p$
\end{condition}

\begin{theorem} \label{Approximate Error Bound}(\textbf{Approximation Error Bound}) Under Condition 3.1-3.3,
\be
\|\bbeta^*(C_u,C_l)-\bbeta^*\|_2\leq \frac{2^k}{(k-1)}\frac{\max\{\alpha,1-\alpha\}}{\min\{\alpha,1-\alpha\}}\frac{\sqrt{\kappa_u}}{\kappa_l}(M_k+c\kappa_0^k)C^{1-k}.
\ee
\end{theorem}
where $C=\min\{C_u,|C_l|\}$ and $c$ is a positive constant only depending on $k$.

Theorem \ref{Approximate Error Bound} indicates that the approximation error vanishes faster if higher moments of $\epsilon^k|\x$ exist. What's more, different from the classical huber estimation, $C_u,|C_l|$ are allowed to diverge so that the approximation cost can be little if $C=\min\{C_u,|C_l|\}$ is chosen properly.

Now focus our attention on the estimation error $\|\hat{\bbeta}-\bbeta^*(C_u,C_l)\|_2$.
Denote by the empirical loss function for the optimization problem (\ref{non-convex regularized problem}):
\be
L_n(\boldsymbol{\beta})=\frac{1}{n}\sum_{i=1}^n\psi_{\alpha}(y_i-\x_i'\boldsymbol{\beta};C_u,C_l).
\ee
 Note that $L_n(\boldsymbol{\beta})$ is convex but not strongly convex due to high dimensionality. To establish the finite sample estimation error bound, we need the so-called restricted strong convexity (RSC). Such condition has been widely discussed in high-dimensional non-convex analysis, for example, see \cite{ANW12} and \cite{LaW15}. This condition imposes a lower bound on the remainder after the first-order Taylor expansion of $L_n(\boldsymbol{\beta})$ and fundamentally requires the curvature not too flat near the optimal minima. For simplicity and convenience in notation, we omit the notation dependence with parameters $C_u,~C_l$ and denote by $\tilde{\bbeta^*}:=\bbeta^*(C_u,C_l)$ for short when no confusion arises.
\begin{definition}\textbf{(Restricted Strong Convexity,~\cite{LaW15})}
Let $\delta L_n(\tilde{\bbeta}^*,\Delta)$ be the remainder after the first-order Taylor expansion around $\tilde{\bbeta}^*$, i.e.,
\be
\delta L_n(\tilde{\bbeta}^*,\Delta)=L_n(\tilde{\bbeta}^*+\Delta)-L_n(\tilde{\bbeta}^*)-\langle\nabla L_n(\tilde{\bbeta}^*),\Delta\rangle.
\ee
The Restricted Strong Convexity condition is defined as follows:
\be\label{restricted strong convexity }
\delta L_n(\tilde{\bbeta}^*,\Delta)\geq
\begin{cases}
\kappa_1 \|\Delta\|_2^2-\tau_1\frac{\log p}{n}\|\Delta\|_1^2, &\forall \|\Delta\|_2\leq 1, \\
\kappa_2 \|\Delta\|_2-\tau_2\sqrt{\frac{\log p}{n}}\|\Delta\|_1, &\forall \|\Delta\|_2 \geq 1.
\end{cases}
\ee
\end{definition}

Since the penalty is non-convex, we introduce one additional condition, the side condition $\|\bbeta\|_1\leq R$ in the problem (\ref{non-convex regularized problem}), where $R$ should be treated carefully. $\|\bbeta^*\|_1\leq R$ is also required so that the true regression vector $\bbeta^*$ is feasible for this problem. This condition originally comes from the requirement of the feasible set in Lasso-type regularized framework. In fact, $L_n(\boldsymbol{\beta})$ and the penalty $P_{\lambda}(\cdot)$ are continuous, with this side condition, the Weierstrass extreme value theorem guarantees that a global minimum for Problem (\ref{non-convex regularized problem}) exists.

Lemma \ref{lemma: restricted strong convexity} and \ref{lemma: restricted strong convexity 2} indicate that the empirical loss $L_n(\boldsymbol{\beta})$ satisfies the Restricted Strong Convexity condition with overwhelming probability. On this basis, we get the following upper bound of the estimation error $\|\hat{\bbeta}-\bbeta^*(C_u,C_l)\|_2$.

\begin{theorem} (\textbf{Estimation Error Bound})\label{estimation error bound}
Suppose Conditions 3.1-3.3 hold and \\ $\|\bbeta\|_1\leq R$,~$n\geq \max\{4R^2 \tau_1^2,\frac{16R^2\max\{\tau_1^2,1\}}{\kappa_1^2}\}\log p$. If the tuning parameter $\lambda$ is chosen to satisfy $\kappa_1\sqrt{\frac{\log p}{n}}\leq \lambda\leq\frac{\kappa_1}{6RL}$, $\max\{C_u,|C_l|\}\leq c/\lambda$ where positive constant $c$ depends on $M_k, \kappa_0$ and $L$,
then there exist positive constants $c_1,~c_2$ such that, with probability at least $1-c_1\exp\{-c_2n\}$,
\be
\parallel\hat{\boldsymbol{\beta}}-\bbeta^*(C_u,C_l)\parallel_2\leq \frac{6\lambda L \sqrt{q}}{4\kappa_1-3\mu}, %~~~~\text{and}~~~
%\parallel\hat{\boldsymbol{\beta}}-\bbeta^*(C_u,C_l)\parallel_1\leq \frac{24\lambda L q}{4\kappa_1-3\mu}.
\ee
\end{theorem}

The result in this theorem can be extended to holds for any vector $\tilde{\bbeta}$ satisfying the first-order necessary condition of problem (\ref{non-convex regularized problem})
\by
\langle\nabla L_n(\tilde{\bbeta})+\nabla P_{\lambda}(\tilde{\bbeta}),\bbeta-\tilde{\bbeta}\rangle\geq 0,\text{for all feasible} ~~\bbeta\in\R^p,
\ey
similar results in the Theorem 1 of \cite{LaW15}. If $\tilde{\bbeta}$ lies in the interior of the feasible set, this necessary condition reduces to the usual formula
$$\nabla L_n(\tilde{\bbeta})+\nabla P_{\lambda}(\tilde{\bbeta})=0.$$
This observation inspires us to use the well-studied first-order algorithms, as we did before in Algorithm \ref{LLA algorithm}.

Finally, Theorems \ref{Approximate Error Bound} and \ref{estimation error bound} together give the nonasymptotic
upper bound of the statistical error $\|\hat{\bbeta}-\bbeta^*\|_2$.
\begin{theorem} (\textbf{Statistical Error Bound})\label{statistical error bound} Under the conditions in Theorems \ref{Approximate Error Bound} and \ref{estimation error bound}, with probability $1-c_1\exp\{-c_2n\}$,
\be
\|\hat{\boldsymbol{\beta}}-\bbeta^*\|_2&\leq& \|\hat{\boldsymbol{\beta}}-\bbeta^*(C_u,C_l)\|_2+\|\bbeta^*(C_u,C_l)-\bbeta^*\|_2\nonumber\\
&\leq& \frac{2^k}{(k-1)}\frac{\max\{\alpha,1-\alpha\}}{\min\{\alpha,1-\alpha\}}\frac{\sqrt{\kappa_u}}{\kappa_l}(M_k+c\kappa_0^k)C^{1-k}+\frac{6\lambda L \sqrt{q}}{4\kappa_1-3\mu}.
\ee
\end{theorem}

By the result in Theorem \ref{statistical error bound}, if $q=q(n)=O(n^{b_1})$ and $\log p=O(n^{b_2})$ for $0<b_1,b_2<1$ and $0<b_1+b_2<1$, then
by the divergent requirement on $\min\{C_u,|C_l|\}$ and the relation,
$\max\{C_u,|C_l|\}\leq c/\lambda$, we can choose $C_u \asymp |C_l|=O\left((\frac{n}{q\log p})^{1/2(k-1)}\right)$ so that the statistical error bound can achieve the optimal convergence rate, $\|\hat{\boldsymbol{\beta}}-\bbeta^*\|_2=O_p(\sqrt{\frac{q\log p}{n}})$.
In practise, the distribution of errors is unknown and $C_u,|C_l|$, along with $\lambda$, can be chosen by  for example, cross validation. We can perform the three-dimensional grid research for the corresponding optimal values, but have to pay the cost of computational time. In the future work, data-driven and tuning-free schemes should be proposed for these parameters.

%Oracle估计量的逼近真值的速度
%估计量的收敛速度
%异方差性的分析
%算法的收敛性
\section{Simulation}
\label{simu}
\setcounter{definition}{0}\setcounter{definition}{0}
\setcounter{equation}{0}\setcounter{lemma}{0}
\setcounter{proposition}{0}\setcounter{theorem}{0}
\setcounter{remark}{0}\setcounter{corollary}{0}
In this section, we assess the finite sample performances of the proposed robust regularized expectile regression. For the choice of the general folded concave penalty function $P_{\lambda}(t)$, here we use the SCAD penalty as an example. The same procedure can be carried out with the MCP penalty or other penalties satisfying  \textbf{Assumption (\ref{assumption on penalty})}. We are mainly engaged in two aspects in this simulation study: one is to investigate the performances of our proposed method in coefficient estimation and model selection; the other is to detect heteroscedasticity by means of regularized expectile regression at different choices of weight levels $\alpha$. For convenience, denote the penalized robust expectile regression with the SCAD penalty by RE-SCAD for short.

We adopt a similar high-dimensional heteroscedastic model from \cite{Wan12} and \cite{ZCZ18}. In this data generation process, firstly, the quasi-covariates $\tilde{\x}=(\tilde{x}_1,\ldots,\tilde{x}_p)'$  are generated from
  the multivariate normal distribution $N_p(\mathbf{0},\Sigma)$ where $\Sigma=(\sigma_{ij})_{p\times p}$, $\sigma_{ij}=0.5^{|i-j|}$ for $i,j=1,\ldots,p$.
   Then we set $x_1=\sqrt{12}\Phi(\tilde{x}_1)$ where $\Phi(\cdot)$ is the cumulative distribution function of the standard normal variable and $\sqrt{12}$ scales $x_1$ to have variance 1. For $i=2,\ldots,p$, $x_i=\tilde{x}_i$. Then the response variable $y$ is generated according to the following sparse heteroscedastic model,
\bq
y=x_6+x_{12}+x_{15}+x_{20}+0.70x_1\epsilon,
\eq
where $\epsilon$ is independent of the covariates $\x$. To investigate how the proposed method performs when the error $\epsilon$ involves various shapes of distribution, we consider the following 4 scenarios:
\bn[(1)]
\item Standard normal distribution $N(0,1)$;
%\item Mixture of normal distribution (MixN):$I\times N(0,1)+(1-I)\times N(0,225)$ where $I$ follows a Bernoulli distribution with probability 0.90; This case is designed where outliers exist.
\item Standard t-distribution with degrees of freedom 3 ($t_3$), where $\ep\epsilon^{2+\delta}<\infty$ exists for $\delta \in (0,1)$.
\item Log-normal distribution (Log-Normal):$\epsilon=\exp\{1.2Z\}$, where $Z$ follows standard normal distribution.
%The corresponding mean and variance of this error distribution is ? and ? respectively.
\item Weibull distribution (Weibull) with shape parameter 0.3 and scale parameter 0.5.
\en

 From this data generation process above, we can see that the true coefficients in the mean part are sparse and only include 4 informative variables. We should also note that $x_1$ plays a role in the error term
 and results in heteroscedasticity. So $x_1$ should be regarded as the significant variable since it plays an essential role in the conditional distribution of $y$ given the covariates. As for the 4 scenarios of error distribution, compared to  the standard normal distribution, standard Student t-distribution $t_3$, Log-normal distribution and Weibull distribution are designed to test the performances of our method when heavy-tailed phenomena exist.

As we point out, $\bbeta^*$ is allowed to change at different weight level, so is in this data generation process due to heteroscedasticity. So we need to perform coefficient estimation and model selection under the identification condition $m_{\alpha}(\epsilon|\x)=0$. In our simulation study, this satisfying $\alpha$ is artificially pre-determined to 0.50 for convenience. Noting that $0.50$-th expectile is just the mean, so the errors in all scenarios are standardized to have mean 0.

 %{\red As we point out, $\bbeta^*$ may change for different weight level. So for some specific  expectile weight level $\alpha$,  we need the identification condition $m_{\alpha}(\epsilon|\x)=0$ to perform coefficient estimation and model selection. In our simulation study, this satisfying $\alpha$ is artificially pre-determined to 0.50 for convenience. With this setting, expectile regression is exactly the classical ordinary least squares regression. At the same time, due to this assumption, the covariates resulting in heteroscedasticity do not get involved in and so are not selected. }

For comparison purpose, in this simulation, we also investigate the performances of the regularized expectile regression in \cite{ZCZ18}, denoted by E-SCAD for short,
\bq\label{E-SCAD}
\underset{\bbeta \in \R^{p}}{\arg\min}\frac{1}{n}\sum_{i=1}^n\phi_{\alpha}(y_i-\x_i'\bbeta)+\sum_{j=1}^pP_{\lambda}(|\beta_j|).
\eq
Similarly as in Algorithm \ref{LLA algorithm}, we can convert problem (\ref{E-SCAD}) to a convex one through  the local linear approximation strategy and use the CVX package to solve it. Besides, for comparison benchmark, we introduce the oracle estimator (\ref{Oracle Estimator Representation}) as the benchmark of estimation accuracy.  Taking heavy-tailed situation into consideration, we define the following robust oracle estimator $\hat{\bbeta}^*=(\hat{\bbeta}_A^{*'},\mathbf{0}_{p-q}')'$,
\be\label{Oracle Estimator Representation}
\hat{\bbeta}_A^*=\underset{\bbeta\in \R^{q}}{\arg\min}\frac{1}{n}\sum_{i=1}^n\psi_{\alpha}(y_i-\z_i'\bbeta;C_u,C_l).
\ee
We aim to show that the introduction of robustness does bring benefits and advantages in model selection consistency at the acceptable cost of coefficient estimation accuracy when the error is heavy-tailed distributed.

 We set the sample size $n=300$ and covariate dimension $p=400$.  Besides $\alpha=0.50$,
 inspired by the results in~\cite{New87}, the positions of weight level near the tail seem to be more effective for testing heteroscedasticity, so we consider two other positions: $\alpha=0.10,0.90$ for this purpose. The combined numerical results are used to test whether our proposed method can detect heteroscedasticity. Given expectile weight level, there are four tuning parameters in total: $a$ and $\lambda$ in SCAD penalty function, and $C_u$ and $C_l$ in robust asymmetric square loss function. We follow the suggestion proposed by \cite{Fan01} and set $a=3.7 $ to reduce the computation burden. For the tuning parameter $\lambda$, $C_u$ and $C_l$,  we adopt the cross validation strategy and perform a three-dimensional grid search on another tuning data set with size $10n$. $\lambda$, $C_u$ and $C_l$ are chosen to minimize the prediction expectile loss error calculated on this tuning data set.

Let $\hat{\bbeta}$ be the coefficient estimates from a given method. We repeat the simulation procedure 100 times and evaluate the
performance in term of the following criteria:
\bi
\item AE:~the average absolute estimation error defined by $\sum_{i=1}^p|\hat{\beta}_j-\beta_j^*|$.
\item SE:~the average square estimation error defined by $\sqrt{\sum_{i=1}^p|\hat{\beta}_j-\beta_j^*|^2}$.
\item Size:~the average number of nonzero regression coefficients $\hat{\beta}_j\neq0$ for $j=1,\ldots,p$. At $\alpha=0.50$, the true size is 4 due to $m_{0.50}(\epsilon|\x)=0$ and at other two positions $\alpha=0.10,~0.90$, given the role of $x_1$, the true size is 5.
\item F:~the frequency that $x_6,x_{12},x_{15},x_{20}$ are selected during the 100 repetitions.
\item F1:~the frequency that $x_1$ is selected during the 100 repetitions.

\ei

\begin{table}[!htp]
	\centering
	\fontsize{8}{10}\selectfont
	\begin{threeparttable}
		\caption{Simulation results under $N(0,1)$ when $n=300$.}
		\label{tab:performance_comparison 1}
		\begin{tabular}{cccccccc}
			\toprule
			&Expectile Level& Method & AE&SE&Size&F&$F_1$\cr	\cmidrule{1-8}
			\multirow{10}{*}{$p=400$}&\multirow{3}{*}{$\alpha=0.10$}&RE-SCAD&0.2449(0.0947)&0.1438(0.0534)&5(0)&100&100\cr
			&&E-SCAD&0.2995(0.1112)&0.1679(0.0578)&5.7700(0.8973)&100&100\cr
			&&Oracle&0.2450(0.0946)&0.1438(0.0534)&5&-&-\cr
			\cmidrule{2-8}
			&\multirow{3}{*}{$\alpha=0.50$}&RE-SCAD&0.2138(0.0810)&0.1211(0.0447)&4.5500(0.7017)&100&14\cr
			&&E-SCAD&0.3104(0.1291)&0.1828(0.0742)&4.7300(0.8475)&100&21\cr
			&&Oracle&0.2057(0.0785)&0.1199(0.0446)&4&-&-\cr
			\cmidrule{2-8}
			&\multirow{3}{*}{$\alpha=0.90$}&RE-SCAD&0.2684(0.0929)&0.1567(0.0508)&5(0)&100&100\cr
			&&E-SCAD&0.3446(0.1418)&0.1836(0.0634)&6.3500(1.52667)&100&100\cr
			&&Oracle&0.2684(0.0929)&0.1567(0.0508)&5&-&-\cr
			\cmidrule{2-8}
		%	\midrule
		%	%	&Expectile Level& Method & AE&SE&Size&F&$F_1$\cr	\cmidrule{2-8}
		%	\multirow{10}{*}{$p=600$}&\multirow{3}{*}{$\alpha=0.10$}&RE-SCAD&1&1&1&1&1\cr
		%	&&E-SCAD&1&1&1&1&1\cr
		%	&&Oracle&1&1&1&1&1\cr
		%	\cmidrule{2-8}
		%	&\multirow{3}{*}{$\alpha=0.50$}&RE-SCAD&1&1&1&1&1\cr
		%	&&E-SCAD&1&1&1&1&1\cr
		%	&&Oracle&1&1&1&1&1\cr
		%	\cmidrule{2-8}
		%	&\multirow{3}{*}{$\alpha=0.90$}&RE-SCAD&1&1&1&1&1\cr
		%	&&E-SCAD&1&1&1&1&1\cr
		%	&&Oracle&1&1&1&1&1\cr
		%	\cmidrule{2-8}
			\bottomrule
		\end{tabular}
		%\small
		%Note:  Standard deviation error in parentheses based on 100 repetitions.
	\end{threeparttable}
\end{table}

\begin{comment}
\begin{table}[!htp]
	\centering
	\fontsize{5.6}{8}\selectfont
	\begin{threeparttable}
		\caption{Simulation results using RE-SCAD, E-SCAD and oracle methods under Mixture Normal when $n=300$.}
		\label{tab:performance_comparison 2}
		\begin{tabular}{cccccccc}
			\toprule
			&Expectile Level& Method & AE&SE&Size&F&$F_1$\cr	\cmidrule{1-8}
			\multirow{10}{*}{$p=400$}&\multirow{3}{*}{$\alpha=0.10$}&RE-SCAD&1&1&1&1&1\cr
			&&E-SCAD&1&1&1&1&1\cr
			&&Oracle&1&1&1&1&1\cr
			\cmidrule{2-8}
			&\multirow{3}{*}{$\alpha=0.50$}&RE-SCAD&1&1&1&1&1\cr
			&&E-SCAD&1&1&1&1&1\cr
			&&Oracle&1&1&1&1&1\cr
			\cmidrule{2-8}
			&\multirow{3}{*}{$\alpha=0.90$}&RE-SCAD&1&1&1&1&1\cr
			&&E-SCAD&1&1&1&1&1\cr
			&&Oracle&1&1&1&1&1\cr
			\cmidrule{2-8}
			\midrule
			%	&Expectile Level& Method & AE&SE&Size&F&$F_1$\cr	\cmidrule{2-8}
			\multirow{10}{*}{$p=600$}&\multirow{3}{*}{$\alpha=0.10$}&RE-SCAD&1&1&1&1&1\cr
			&&E-SCAD&1&1&1&1&1\cr
			&&Oracle&1&1&1&1&1\cr
			\cmidrule{2-8}
			&\multirow{3}{*}{$\alpha=0.50$}&RE-SCAD&1&1&1&1&1\cr
			&&E-SCAD&1&1&1&1&1\cr
			&&Oracle&1&1&1&1&1\cr
			\cmidrule{2-8}
			&\multirow{3}{*}{$\alpha=0.90$}&RE-SCAD&1&1&1&1&1\cr
			&&E-SCAD&1&1&1&1&1\cr
			&&Oracle&1&1&1&1&1\cr
			\cmidrule{2-8}
			\bottomrule
		\end{tabular}
		%\small
		%Note:  Standard deviation error in parentheses based on 100 repetitions.
	\end{threeparttable}
\end{table}
\end{comment}
\begin{table}[!htp]
	\centering
	\fontsize{8}{10}\selectfont
	\begin{threeparttable}
		\caption{Simulation results under $t_3$ when $n=300$.}
		\label{tab:performance_comparison 2}
		\begin{tabular}{cccccccc}
			\toprule
			&Expectile Level& Method & AE&SE&Size&F&$F_1$\cr	\cmidrule{1-8}
			\multirow{10}{*}{$p=400$}&\multirow{3}{*}{$\alpha=0.10$}&RE-SCAD&0.6427(0.6734)&0.3965(0.4170)&4.8400(0.5265)&90&100\cr
			&&E-SCAD&2.0011(2.5921)&0.8302(0.6566)&10.9100(6.6835)&89&100\cr
			&&Oracle&0.4037(0.1669)&0.2381(0.0965)&5&-&-\cr
			\cmidrule{2-8}
			&\multirow{3}{*}{$\alpha=0.50$}&RE-SCAD&0.2533(0.1788)&0.1518(0.1211)&4.0200(0.1407)&100&2\cr
			&&E-SCAD&1.2436(0.8276)&10.4386(0.1940)&15.5900(8.4305)&100&36\cr
			&&Oracle&0.2311(0.0873)&0.1353(0.0506)&4&-&-\cr
			\cmidrule{2-8}
			&\multirow{3}{*}{$\alpha=0.90$}&RE-SCAD&0.4312(0.1990)&0.2545(0.1251)&5.0400(0.1969)&100&100\cr
			&&E-SCAD&2.0352(3.0028)&0.8284(0.6666)&11.7000(7.4678)&92&100\cr
			&&Oracle&0.4250(0.1856)&0.2481(0.1008)&5&-&-\cr
			\cmidrule{2-8}
			%\midrule
			%	&Expectile Level& Method & AE&SE&Size&F&$F_1$\cr	\cmidrule{2-8}
	%		\multirow{10}{*}{$p=600$}&\multirow{3}{*}{$\alpha=0.10$}&RE-SCAD&1&1&1&1&1\cr
	%		&&E-SCAD&1&1&1&1&1\cr
	%		&&Oracle&1&1&1&1&1\cr
	%		\cmidrule{2-8}
	%		&\multirow{3}{*}{$\alpha=0.50$}&RE-SCAD&1&1&1&1&1\cr
	%		&&E-SCAD&1&1&1&1&1\cr
	%		&&Oracle&1&1&1&1&1\cr
	%		\cmidrule{2-8}
	%		&\multirow{3}{*}{$\alpha=0.90$}&RE-SCAD&1&1&1&1&1\cr
	%		&&E-SCAD&1&1&1&1&1\cr
	%		&&Oracle&1&1&1&1&1\cr
	%		\cmidrule{2-8}
			\bottomrule
		\end{tabular}
		%\small
		%Note:  Standard deviation error in parentheses based on 100 repetitions.
	\end{threeparttable}
\end{table}

\begin{table}[!htp]
	\centering
	\fontsize{8}{10}\selectfont
	\begin{threeparttable}
		\caption{Simulation results  under $Log-normal$ when $n=300$.}
		\label{tab:performance_comparison 3}
		\begin{tabular}{cccccccc}
			\toprule
			&Expectile Level& Method & AE&SE&Size&F&$F_1$\cr	\cmidrule{1-8}
			\multirow{10}{*}{$p=400$}&\multirow{3}{*}{$\alpha=0.10$}&RE-SCAD
			&0.1113(0.0420)&0.0650(0.0.0240)&5(0)&100&100\cr
			&&E-SCAD&0.9492(4.6602)&0.2794(0.7496)&8.100(9.9186)&98&100\cr
			&&Oracle&0.1113(0.0420)&0.0650(0.0.0240)&5&-&-\cr
			\cmidrule{2-8}
			&\multirow{3}{*}{$\alpha=0.50$}&RE-SCAD&0.2212(0.0835)&0.1280(0.0466)&5.3200(0.5840)&100&100\cr
			&&E-SCAD&9.2385(24.4768)&1.9726(2.8154)&26.7900(22.5685)&77&56\cr
			&&Oracle&0.2153(0.0822)&0.1272(0.0466)&4&-&-\cr
			\cmidrule{2-8}
			&\multirow{3}{*}{$\alpha=0.90$}&RE-SCAD&0.5365(0.2815)&0.2957(0.1524)&6.1000(1.6484)&100&68\cr
			&&E-SCAD&13.5400(16.5492)&3.6973(2.6933)&21.1400(11.7774)&1&100\cr
			&&Oracle&0.4672(0.1796)&0.2763(0.1000)&5&-&-\cr
			\cmidrule{2-8}
			% \midrule
			%	&Expectile Level& Method & AE&SE&Size&F&$F_1$\cr	\cmidrule{2-8}
		%	\multirow{10}{*}{$p=600$}&\multirow{3}{*}{$\alpha=0.10$}&RE-SCAD&1&1&1&1&1\cr
		%	&&E-SCAD&1&1&1&1&1\cr
		%	&&Oracle&1&1&1&1&1\cr
		%	\cmidrule{2-8}
		%	&\multirow{3}{*}{$\alpha=0.50$}&RE-SCAD&1&1&1&1&1\cr
		%	&&E-SCAD&1&1&1&1&1\cr
		%	&&Oracle&1&1&1&1&1\cr
		%	\cmidrule{2-8}
		%	&\multirow{3}{*}{$\alpha=0.90$}&RE-SCAD&1&1&1&1&1\cr
		%	&&E-SCAD&1&1&1&1&1\cr
		%	&&Oracle&1&1&1&1&1\cr
		%	\cmidrule{2-8}
			\bottomrule
		\end{tabular}
		%\small
		%Note:  Standard deviation error in parentheses based on 100 repetitions.
	\end{threeparttable}
\end{table}

\begin{table}[!htp]
	\centering
	\fontsize{8}{10}\selectfont
	\begin{threeparttable}
		\caption{Simulation results under Weibull when $n=300$.}
		\label{tab:performance_comparison 4}
		\begin{tabular}{cccccccc}
			\toprule
			&Expectile Level& Method & AE&SE&Size&F&$F_1$\cr	\cmidrule{1-8}
			\multirow{10}{*}{$p=400$}&\multirow{3}{*}{$\alpha=0.10$}&RE-SCAD&0.0371(0.0171)&0.0216(0.0097)&5(0)&100&100\cr
			&&E-SCAD&NA&NA&NA&NA&NA\cr
			&&Oracle&0.0371(0.0171)&0.0216(0.0097)&5&-&-\cr
			\cmidrule{2-8}
			&\multirow{3}{*}{$\alpha=0.50$}&RE-SCAD&0.1708(0.0697)&0.0996(0.0380)&5.3400(0.5724)&100&100\cr
			&&E-SCAD&NA&NA&NA&NA&NA\cr
			&&Oracle&0.1651(0.0670)&0.0980(0.0381)&4&-&-\cr
			\cmidrule{2-8}
			&\multirow{3}{*}{$\alpha=0.90$}&RE-SCAD&0.8253(0.5898)&0.4802(0.3247)&6.1400(1.1636)&95&100\cr
			&&E-SCAD&NA&NA&NA&NA&NA\cr
			&&Oracle&0.5936(0.2249)&0.3551(0.1252)&5&-&-\cr
			\cmidrule{2-8}
		%	\midrule
			%	&Expectile Level& Method & AE&SE&Size&F&$F_1$\cr	\cmidrule{2-8}
		%	\multirow{10}{*}{$p=600$}&\multirow{3}{*}{$\alpha=0.10$}&RE-SCAD&1&1&1&1&1\cr
		%	&&E-SCAD&1&1&1&1&1\cr
		%	&&Oracle&1&1&1&1&1\cr
		%	\cmidrule{2-8}
		%	&\multirow{3}{*}{$\alpha=0.50$}&RE-SCAD&1&1&1&1&1\cr
		%	&&E-SCAD&1&1&1&1&1\cr
		%	&&Oracle&1&1&1&1&1\cr
		%	\cmidrule{2-8}
		%	&\multirow{3}{*}{$\alpha=0.90$}&RE-SCAD&1&1&1&1&1\cr
		%	&&E-SCAD&1&1&1&1&1\cr
		%	&&Oracle&1&1&1&1&1\cr
		%	\cmidrule{2-8}
			\bottomrule
		\end{tabular}
		\small
		Note:  'NA' stands for Not Available;
		%Standard deviation error in parentheses based on 100 repetitions.
	\end{threeparttable}
\end{table}

When the error follows normal distribution, RE-SCAD and E-SCAD perform very similarly in both estimation accuracy and model recognition, see Table \ref{tab:performance_comparison 1}. But things get changed when heavy-tailed phenomena arise. RE-SCAD always performs well, showing robustness to various distributions, while E-SCAD gets worse (Table \ref{tab:performance_comparison 2},\ref{tab:performance_comparison 3}), even fails (Table \ref{tab:performance_comparison 4}). In these situations, the classical E-SCAD without robust modification losses estimation efficiency and even fails to pick up the informative variables in the mean part (see 'F' column), not to mention poor performances in detecting heteroscedasticity (see 'F1' column). An interesting founding is that when the error is asymmetric and heavy-tailed like the Log-Normal or Weibull distribution, robust expectile regressions at $\alpha=0.10$ and $\alpha=0.90$ appear different performances in coefficient estimation and model recognition, indicating great potential of expectile regression when asymmetry exists.
From all the discussion above, we can come to an conclusion that our method shares satisfactory performances for modelling high-dimensional heavy-tailed and heterogeneous data.

\begin{comment}

\subsection{An real Example}
Gene expression data for predicting log transformed riboflavin (Vitamin B2) production rate in Bacillus subtilis. (\cite{BKM14}). It contains 71 observations and 4088 features (gene expressions). This dataset is available in R package $hdi$ via $data(riboflavin)$. For the task, only 1000 features with the largest variance were used. This step can be completed by some kind of sure screening methods ???
\end{comment}
\section{Appendix}
\label{Appe}
\setcounter{definition}{0}\setcounter{definition}{0}
\setcounter{equation}{0}\setcounter{lemma}{0}
\setcounter{proposition}{0}\setcounter{theorem}{0}
\setcounter{remark}{0}\setcounter{corollary}{0}
\begin{lemma}\label{phi strong convex}
The loss function $\phi(\cdot)$ defined in (\ref{phi}) is continuous differentiable. Moreover, for any $r,r_0\in \R$, we have
\bq
\min\{\alpha,1-\alpha\}\cdot(r-r_0)^2\leq \phi_{\alpha}(r)-\phi_{\alpha}(r_0)-\phi'_{\alpha}(r_0)\cdot(r-r_0)\leq \max\{\alpha,1-\alpha\}\cdot(r-r_0)^2.
\eq
%where $\underline{\alpha}=\min\{\alpha,1-\alpha\}$ and $\bar{\alpha}=\max\{\alpha,1-\alpha\}$.
%It follows that $\phi(\cdot)$ is strongly convex.
\end{lemma}
\proof Details can be found in \cite{GaZ16}.

\textbf{Proof of Theorem \ref{Approximate Error Bound}}
 For simplicity and convenience in notation, we omit the notation dependence with the pre-determined parameters $C_u,~C_l$ and denote by $\tilde{\bbeta^*}:=\bbeta^*(C_u,C_l)$ for short.

By equation (\ref{theoretical optimal solution condition}) and Lemma \ref{phi strong convex},
\be
\ep[\phi_{\alpha}(y-\x'\tilde{\bbeta}^*)-\phi_{\alpha}(y-\x'\bbeta^*)]&\geq& \min\{\alpha,1-\alpha\}\ep[\x'\tilde{\bbeta}^*-\x'\bbeta^*]^2\nonumber\\
&=&\min\{\alpha,1-\alpha\}(\tilde{\bbeta}^*-\bbeta^*)'\ep(\x\x')(\tilde{\bbeta}^*-\bbeta^*)\nonumber\\
&\geq& \min\{\alpha,1-\alpha\}\kappa_l\parallel\tilde{\bbeta}^*-\bbeta^*\parallel_2^2
\ee
where the last inequality follows by Condition \ref{C(2)}.

On the other hand, $\tilde{\bbeta}^*=\underset{\boldsymbol{\beta} \in \R^{p}}{\arg\min}~\ep\psi_{\alpha}(y-\x'\bbeta;C_u,C_l)$.  Then
\by
&&\ep[\phi_{\alpha}(y-\x'\tilde{\bbeta}^*)-\phi_{\alpha}(y-\x'\bbeta^*)]\nonumber\\
&=&\ep[\phi_{\alpha}(y-\x'\tilde{\bbeta}^*)-\psi_{\alpha}(y-\x'\tilde{\bbeta}^*)]
+\ep[\psi_{\alpha}(y-\x'\tilde{\bbeta}^*)-\psi_{\alpha}(y-\x'\bbeta^*)]
-\ep[\phi_{\alpha}(y-\x'\bbeta^*)-\psi_{\alpha}(y-\x'\bbeta^*)]\nonumber\\
&\leq& \ep[g_{\alpha}(y-\x'\tilde{\bbeta}^*)-g_{\alpha}(y-\x'\bbeta^*)],
\ey
where $g_{\alpha}(\cdot)$ is defined as follows
\be
g_{\alpha}(r)=\phi_{\alpha}(r)-\psi_{\alpha}(r)=\alpha (r-C_u)^2\I(r\geq C_u)+(1-\alpha)(r-C_l)^2\I(r\leq C_l).
\ee
Note $g_{\alpha}(r)$ is continuous and differentiable and
\bqs
g'_{\alpha}(r)=2\alpha (r-C_u)\I(r\geq C_u)+2(1-\alpha)(r-C_l)\I(r\leq C_l).
\eqs
So by the mean value theorem, there exists some $\tilde{\bbeta}$ on the line segment between $\tilde{\bbeta}^*$ and $\bbeta^*$ such that
\by
\left|\ep[g_{\alpha}(y-\x'\tilde{\bbeta}^*)-g_{\alpha}(y-\x'\bbeta^*)]\right|&=& \ep\left[|g_{\alpha}'(y-\x\tilde{\bbeta})|\times|\x'(\tilde{\bbeta}^*-\bbeta^*)|\right]\nonumber\\
&\leq& 2\max\{\alpha,1-\alpha\}\ep\left[(\tilde{r}-C)\I(\tilde{r}\geq C)\times|\x'(\tilde{\bbeta}^*-\bbeta^*)|\right]
\ey
where $\tilde{r}=|y-\x'\tilde{\bbeta}|$ and $C=\min\{C_u,|C_l|\}$.

Denote by $P_{\epsilon}$ the conditional distribution of $\epsilon$ on $\x$ and $\ep_{\epsilon}$ the corresponding conditional expectation, we have
\be
\ep_{\epsilon}(\tilde{r}-C)\I(\tilde{r}\geq C)&=&\int_{0}^{\infty}P_{\epsilon}(\tilde{r}\I(\tilde{r}\geq C)>t)dt-CP_{\epsilon}(\tilde{r}\geq C)\nonumber\\
&=& \int_{0}^{\infty}P_{\epsilon}(\tilde{r}>t,~\tilde{r}>C)dt-CP_{\epsilon}(\tilde{r}\geq C)\nonumber\\
&=& \int_{C}^{\infty}P_{\epsilon}(\tilde{r}>t)dt+\int_{0}^{C}P_{\epsilon}(\tilde{r}>C)dt-CP_{\epsilon}(\tilde{r}\geq C)\nonumber\\
&\leq& \int_{C}^{\infty}\frac{\ep_{\epsilon}[\tilde{r}]^k}{t^k}dt=\frac{1}{k-1}C^{1-k}\ep_{\epsilon}[\tilde{r}]^k
\ee
where the second to last inequality is obtained by Markov inequality.

Therefore, $\ep[\phi_{\alpha}(y-\x'\tilde{\bbeta}^*)-\phi_{\alpha}(y-\x'\bbeta^*)]$ is further bounded by
\be
&&\frac{2}{k-1}C^{1-k}\ep\left[|y-\x'\tilde{\bbeta}|^k|\x'(\tilde{\bbeta}^*-\bbeta^*)|\right]\nonumber\\
&=&\frac{2}{k-1}C^{1-k}\ep\left[|\epsilon+\x'(\bbeta^*-\tilde{\bbeta})|^k|\x'(\tilde{\bbeta}^*-\bbeta^*)|\right]\nonumber\\
&\leq& \frac{2}{k-1}(\frac{C}{2})^{1-k}\ep\left[(|\epsilon|^k+|\x'(\bbeta^*-\tilde{\bbeta})|^k)|\x'(\tilde{\bbeta}^*-\bbeta^*)|\right]\nonumber\\
&=&\frac{2}{k-1}(\frac{C}{2})^{1-k}\left\{\ep\left[|\epsilon|^k|\x'(\tilde{\bbeta}^*-\bbeta^*)|\right]+\ep\left[|\x'(\bbeta^*-\tilde{\bbeta})|^k|\x'(\tilde{\bbeta}^*-\bbeta^*)|\right]\right\}
\ee
As for the first term, by Condition \ref{C(1)} and \ref{C(2)},
\be
\ep\left[|\epsilon|^k|\x'(\tilde{\bbeta}^*-\bbeta^*)|\right]&=&\ep\left[\ep(|\epsilon|^k|\x)|\x'(\tilde{\bbeta}^*-\bbeta^*)|\right]\nonumber\\
&\leq& \left[\ep(\ep(|\epsilon|^k|\x))^2\right]^{1/2}\left[\ep(|\x'(\tilde{\bbeta}^*-\bbeta^*)|^2)\right]^{1/2}\nonumber\\
&\leq&M_k\sqrt{\kappa_u}\parallel\tilde{\bbeta}^*-\bbeta^*\parallel_2
\ee
As for the second term,
\be
\ep\left[|\x'(\bbeta^*-\tilde{\bbeta})|^k|\x'(\tilde{\bbeta}^*-\bbeta^*)|\right]&\leq& [\ep(|\x'(\bbeta^*-\tilde{\bbeta})|^{2k})]^{1/2}[\ep(|\x'(\tilde{\bbeta}^*-\bbeta^*)|)^2]^{1/2}\nonumber\\
&\leq& c\kappa_0^k\sqrt{\kappa_u}\parallel\tilde{\bbeta}^*-\bbeta^*\parallel_2
\ee
since $\x'(\bbeta^*-\tilde{\bbeta})$ is sub-Gaussian distributed by Condition \ref{C(3)}, hence its $2k$th moment is bounded by $c^2\kappa_0^{2k}$ for a universal positive constant $c$ depending on $k$ only.

Combined with these results, we have
\be
\|\bbeta^*(C_u,C_l)-\bbeta^*\|_2\leq \frac{2^k}{(k-1)}\frac{\max\{\alpha,1-\alpha\}}{\min\{\alpha,1-\alpha\}}\frac{\sqrt{\kappa_u}}{\kappa_l}(M_k+c\kappa_0^k)C^{1-k}.
\ee
So far, the proof has been completed. \qed

\begin{lemma}\label{Bernstein inequality}
Let $\xi_1,\xi_2,\ldots,\xi_n$ be independent real valued random variables. Assume that there exists some positive numbers $\nu$ and $c$ such that
\by
 \sum_{i=1}^n\ep[\xi_i^2]\leq \nu,
 \ey
 and for all integers $k\geq 3$
 \by
 \sum_{i=1}^n\ep[(\xi_i)_+^k]\leq \frac{k!}{2}\nu c^{k-2}.
 \ey
Let $S_n=\sum_{i=1}^n(\xi_i-\ep[\xi_i])$, then for every positive $x$,
\by
\P\left(S_n\geq\sqrt{2\nu x}+cx\right)\leq \exp(-x).
\ey
\end{lemma}
\proof Details can be found in Proposition 2.9 of \cite{Mas07}.

\begin{lemma}\label{lemma: restricted strong convexity}
Under Condition 3.1-3.3, there exist universal positive constants $\kappa_1',\kappa_2', c_1', c_2'$ such that
 with probability at least $1-c_1'\exp{(-c_2'n)}$,
\be
\delta L_n(\bbeta,\Delta)\geq \kappa_1'\|\Delta\|_2\left[\|\Delta\|_2-\kappa_2'\sqrt{\frac{\log p}{n}}\|\Delta\|_1\right]
\ee
for all $\|\bbeta\|_2\leq 4\rho_2$, $\|\Delta\|_2\leq 8\rho_2$ where $\rho_2$ is a sufficiently large constant depending on $R$ and $C=\max\{C_u,|C_l|\}\geq c_u\rho_2^{-1}$ where $c_u$ is a positive constant depending on $M_k,\kappa_l,\kappa_u$ and $\kappa_0$.
\end{lemma}
\proof  Define the set
 $ A=\left\{(\bbeta,\Delta):\|\bbeta\|_2\leq 4\rho_2, ~\|\Delta\|_2\leq 8\rho_2\right\}$, then we can show that for any $(\bbeta,\Delta)\in A$,
 \be\label{contraction inequality}
 \delta L_n(\bbeta,\Delta)\geq \frac{\min\{\alpha,1-\alpha\}}{n}\sum_{i=1}^n \varphi_{\tau \|\Delta\|_2}(\x_i'\Delta\I(|y_i-\x_i'\bbeta|<T).
 \ee
for some proper chosen $T$ and $\tau$ satisfying $T+8\tau\rho_2\leq \min\{C_u,|C_l|\}$, where the threshold function
\be
\varphi_t(u)=u^2\I(|u|\leq t^2/2)+(t-|u|)^2\I(t/2\leq |u|\leq t^2).
\ee

To show (\ref{contraction inequality}), if $|y_i-\x_i'\bbeta|>T$ or $|\x_i'\Delta|>\tau \|\Delta\|_2$, the right hand side of (\ref{contraction inequality}) is 0. So by convexity of the robust loss function (\ref{huber phi}),
(\ref{contraction inequality}) holds trivially. If $|y_i-\x_i'\bbeta|\leq T$ and $|\x_i'\Delta|\leq \tau \|\Delta\|_2$, then by Lemma \ref{phi strong convex},
\be
 \delta L_n(\bbeta,\Delta)&\geq &\frac{\min\{\alpha,1-\alpha\}}{n}\sum_{i=1}^n (\x_i'\Delta)^2\nonumber\\
 &\geq &\frac{\min\{\alpha,1-\alpha\}}{n}\sum_{i=1}^n \varphi_{\tau \|\Delta\|_2}(\x_i'\Delta\I(|y_i-\x_i'\bbeta|<T).
\ee

Based on this inequality for $\delta L_n(\bbeta,\Delta)$, we follow the similar proof procedure of Lemma 2 in \cite{FLW17} and obtain the wanted results.
\begin{lemma}\label{lemma: restricted strong convexity 2}
Suppose $L_n(\boldsymbol{\beta})$ is convex and $\tilde{\bbeta}^*$,~$\tilde{\bbeta}^*+\Delta$ lies in the feasible set so that $\|\Delta\|_1\leq 2R$. If the Restricted Strong Convexity condition holds for $\|\Delta\|_2\leq 1$ and $n\geq 4R^2 \tau_1^2\log p$, then
\be
\delta L_n(\tilde{\bbeta}^*,\Delta)\geq \kappa_1 \|\Delta\|_2-\sqrt{\frac{\log p}{n}}\|\Delta\|_1, ~~~~\forall \|\Delta\|_2 \geq 1
\ee
\end{lemma}
\proof Details can be found in Lemma 8 of \cite{LaW15}.

\textbf{Proof of Theorem \ref{estimation error bound}}: If we can prove the following two claims, then by the Theorem 1 of \cite{LaW15}, this theorem holds:
\bi
\item Claim I: $\|\nabla L_n(\tilde{\bbeta}^*)\|_{\infty}\leq \lambda L/4$ with overwhelming probability;
\item Claim II: the empirical loss $L_n(\boldsymbol{\beta})$ satisfies the Restricted Strong Convexity condition.
\ei

For Claim I, we use the Bernstein inequality (Lemma \ref{Bernstein inequality}) and the union bound to establish this result. Through straight calculation,
\be
\nabla L_n(\tilde{\bbeta}^*)=-\frac{1}{n}\sum_{i=1}^n\psi'_{\alpha}(y_i-\x_i'\tilde{\bbeta}^*;C_u,C_l)\x_i,
\ee
where $\psi'_{\alpha}(\cdot)$ is defined in equation \ref{huber phi derivative}.

Note that $\left|\psi'_{\alpha}(r;C_u,C_l)\right|\leq 2\max\{\alpha,1-\alpha\}|r|$ so that for $j=1,\ldots,p$,
\be
\left|\psi'_{\alpha}(y_i-\x_i'\tilde{\bbeta}^*;C_u,C_l)\x_{ij}\right|\leq 2\max\{\alpha,1-\alpha\}\left|y_i-\x_i'\boldsymbol{\beta}\right|\left|\x_{ij}\right|.
\ee
Then
\be
\ep[\left|\psi'_{\alpha}(y_i-\x_i'\tilde{\bbeta}^*;C_u,C_l)\x_{ij}\right|^2]&\leq& 4 \ep[\left|y_i-\x_i'\tilde{\bbeta}^*\right|^2\left|\x_{ij}\right|^2]\nonumber \\
&\leq& 8\ep\left[\left(\epsilon_i^2+\left|\x_i'(\tilde{\bbeta}^*-\bbeta^*\right|^2\right)\left|\x_{ij}\right|^2\right]\nonumber\\
&=&8\ep\left[\ep(\epsilon_i^2|\x_i)\x_{ij}^2+\left|\x_i'(\tilde{\bbeta}^*-\bbeta^*\right|^2\left|\x_{ij}\right|^2\right]\nonumber\\
&\leq& \nu.
\ee
where the last inequality above follows a similar procedure as in the proof of Theorem \ref{Approximate Error Bound} and $\nu$ is a constant depending on $\kappa_0$ and $M_2$.

 Denote by $C=2\max\{\alpha,1-\alpha\}\times\max\{C_u,|C_l|\}$. Let $A=\frac{|\psi'_{\alpha}(y_i-\x_i'\tilde{\bbeta}^*;C_u,C_l)|}{C}$, then $|A|\leq 1$. By the Theorem 2.7 of \cite{Riv12}, $A\x_{ij}$ is also sub-gaussian with the same parameter $\kappa_0$. Then for any $k\geq 3$, by the Proposition 3.2 of \cite{Riv12}, we have
\be
\ep\left[\left|\psi'_{\alpha}(y_i-\x_i'\tilde{\bbeta}^*;C_u,C_l)\x_{ij}\right|^k\right]\leq C^k \ep[|A\x_{ij}|^k]
\leq \frac{k!}{2}(CB)^{k-2}\nu,
\ee
where $B$ is a constant depending on $\kappa_0$.

Meanwhile by the definition of $\tilde{\bbeta}^*$, $\ep[\psi'_{\alpha}(y_i-\x_i'\tilde{\bbeta}^*;C_u,C_l)\x_{ij}]=0$ for $j=1,\ldots,p$.
Then by the Bernstein inequality from Lemma \ref{Bernstein inequality}, we have for $j=1,\ldots,p$
\be
\P\left( \left|-\frac{1}{n}\sum_{i=1}^n\psi'_{\alpha}(y_i-\x_i'\tilde{\bbeta}^*;C_u,C_l)\x_{ij}\right|\geq\sqrt{2\frac{\nu}{n}t}+\frac{CB}{n}t\right)\leq\exp(-t).
\ee
Chose $t=\frac{n\lambda^2L^2}{128\nu}$ and $C\leq \frac{16\nu}{B\lambda L}$, we have $\frac{CBt}{n}\leq \sqrt{\frac{2\nu t}{n}}$ and therefore,
\be
\P\left( \left|-\frac{1}{n}\sum_{i=1}^n\psi'_{\alpha}(y_i-\x_i'\tilde{\bbeta}^*;C_u,C_l)\x_{ij}\right|\geq \frac{\lambda L}{4}\right)\leq\exp\left(-\frac{n\lambda^2L^2}{128\nu}\right).
\ee

 Through the union bound argument, we have
\be
\P\left( \left\|-\frac{1}{n}\sum_{i=1}^n\psi'_{\alpha}(y_i-\x_i'\tilde{\bbeta}^*;C_u,C_l)\x_{i}\right\|_{\infty}\geq \frac{\lambda L}{4}\right)&\leq& p\exp\left(-\frac{n\lambda^2L^2}{128\nu}\right)\nonumber\\
&=&\exp\left(-\frac{n\lambda^2L^2}{128\nu}+\log p\right).
\ee
Chose $\lambda=\kappa_{\lambda}\sqrt{\frac{\log p}{n}}$ and $\frac{\kappa^2_{\lambda}L^2}{128\nu}-1>0$ such that
$\exp\left(-\frac{n\lambda^2L^2}{128\nu}+\log p\right)\leq \exp \left(-cn\right)$ for some positive $c=\frac{\kappa^2_{\lambda}L^2}{128\nu}-1>0$.

For Claim II, by Lemma \ref{lemma: restricted strong convexity}, for $\|\Delta\|_2\leq 8\rho_2$, with probability at least $1-c_1'\exp{(-c_2'n)}$,
\be
\delta L_n(\bbeta,\Delta)\geq \kappa_1'\|\Delta\|_2^2-\kappa_1'\kappa_2\sqrt{\frac{\log p}{n}}\|\Delta\|_1\|\Delta\|_2.
\ee
Using the fact that $ab\leq (a^2+b^2)/2$, we obtain that
\be
\delta L_n(\bbeta,\Delta)&\geq& \kappa_1'\|\Delta\|_2^2-\left(\frac{1}{2}\kappa_1'\|\Delta\|_2^2+\frac{1}{2}\kappa_1'\kappa_2^2\frac{\log p}{n}\|\Delta\|_1^2\right)\nonumber\\
&=& \kappa_1\|\Delta\|_2^2-\tau_1\frac{\log p}{n}\|\Delta\|_1^2.
\ee
with $\kappa_1=\frac{1}{2}\kappa_1',~~\tau_1=\frac{1}{2}\kappa_1'\kappa_2^2$.

 Without loss of generality, we assume $\rho_2\geq 1/8$. So we have proved that the first scenario of the Restricted Strong Convexity, i.e., the Restricted Strong Convexity condition holds for the empirical loss $L_n(\bbeta)$ when $\|\Delta\|_2\leq 1$. By Lemma \ref{lemma: restricted strong convexity 2}, when $n\geq 4R^2 \tau_1^2\log p$, the whole Restricted Strong Convexity condition holds.

 Then based on the two claims above, by probability union bound argument, there exist positive constant $c_1,c_2$ such that with probability at least $1-c_1\exp\{-c_2n\}$, the statistical error bound holds.

\end{document}